\documentclass[a4paper]{amsart}

\usepackage{amssymb}    % AMS fonts
\usepackage{array}
\usepackage{enumitem}
\usepackage{url}
\usepackage{verbatim}

\theoremstyle{plain}
\newtheorem{theorem}{Theorem}[section]
\newtheorem{lemma}[theorem]{Lemma}
\newtheorem{definition}[theorem]{Definition}

\newtheorem{corollary}[theorem]{Corollary}
\theoremstyle{remark}
\newtheorem{example}[theorem]{Example}
\newtheorem{remark}[theorem]{Remark}
\newtheorem{convention}[theorem]{Convention}
\newtheorem*{acknowledgment}{Acknowledgment}

\numberwithin{equation}{section}

\newcommand{\seclabel}[1]{\label{sec:#1}} % section
\newcommand{\thmlabel}[1]{\label{thm:#1}} % theorem
\newcommand{\lemlabel}[1]{\label{lem:#1}} % lemma
\newcommand{\corlabel}[1]{\label{cor:#1}} % corollary
\newcommand{\deflabel}[1]{\label{def:#1}} % definition
\newcommand{\exmlabel}[1]{\label{exm:#1}} % example
\newcommand{\tablabel}[1]{\label{tab:#1}} % table
\newcommand{\cnvlabel}[1]{\label{cnv:#1}} % convention
\newcommand{\eqnlabel}[1]{\label{eqn:#1}} % equation

\newcommand{\secref}[1]{\ref{sec:#1}} % section
\newcommand{\thmref}[1]{\ref{thm:#1}} % theorem
\newcommand{\lemref}[1]{\ref{lem:#1}} % lemma
\newcommand{\corref}[1]{\ref{cor:#1}} % corollary
\newcommand{\defref}[1]{\ref{def:#1}} % definition
\newcommand{\exmref}[1]{\ref{exm:#1}} % example
\newcommand{\tabref}[1]{\ref{tab:#1}} % table
\newcommand{\cnvref}[1]{\ref{cnv:#1}} % convention
\newcommand{\eqnref}[1]{\eqref{eqn:#1}} % parenthesized eqn ref

\newcommand{\ldiv}{\backslash} % left division
\newcommand{\rdiv}{/} % right division
\newcommand{\inv}{^{-1}}

\DeclareMathOperator{\head}{head}
\DeclareMathOperator{\tail}{tail}

\newcommand{\Qgp}{\textbf{Q}}
\newcommand{\LL}{\textbf{L$\Lambda$}}
\newcommand{\RL}{\textbf{R$\Lambda$}}
\newcommand{\Loop}{\textbf{$\Lambda$}}
\newcommand{\pQ}{\textbf{pQ}}
\newcommand{\Qi}{\textbf{Qi}}
\newcommand{\eQ}{\textbf{eQ}}
\newcommand{\Qe}{\textbf{Qe}}
\newcommand{\Loope}{\textbf{Q1}}
\newcommand{\RPQ}{\textbf{RPQ}}
\newcommand{\RPLL}{\textbf{RPL$\Lambda$}}
\newcommand{\RPRL}{\textbf{RPR$\Lambda$}}
\newcommand{\RPL}{\textbf{RP$\Lambda$}}
\newcommand{\RPV}{\textbf{RP$\mathbb{V}$}}
\newcommand{\RPpQ}{\textbf{RPpQ}}
\newcommand{\RPQi}{\textbf{RPQi}}
\newcommand{\RPeQ}{\textbf{RPeQ}}
\newcommand{\RPQe}{\textbf{RPQe}}
\newcommand{\RPLe}{\textbf{RPQ1}}

\newcommand{\Aut}{\mathrm{Aut}} % automorphism group
\newcommand{\End}{\mathrm{End}} % endomorphism semigroup
\newcommand{\Hom}{\mathrm{Hom}} % morphisms

\newcommand{\cR}{\mathcal{R}}

\newcommand{\bN}{\mathbb{N}}

\newcommand{\eqA}{\textsc{A}}
\newcommand{\eqQ}{\textsc{Q}}
\newcommand{\eqB}{\textsc{B}}
\newcommand{\eqV}{\textsc{V}}
\newcommand{\eqhV}{\textsc{\^{V}}}

\newcommand{\Sub}{\operatorname{Sub}}
\newcommand{\Con}{\operatorname{Con}}
\newcommand{\Eq}{\operatorname{Eq}}
\newcommand{\Var}{\operatorname{Var}}
\newcommand{\Free}{\operatorname{Free}}

\title[Right Product Quasigroups and Loops]
{Right Product Quasigroups and Loops}

\author[M.~K.~Kinyon]{Michael~K.~Kinyon}
\address{Department of Mathematics \\
University of Denver \\
2360 S Gaylord St
\\ Denver, Colorado 80208 USA}
\email{\url{mkinyon@math.du.edu}}
\urladdr{\url{http://www.math.du.edu/~mkinyon}}
\author [A. Krape\v{z}]{Aleksandar Krape\v{z}$^*$}
\thanks{$^*$Supported by the Ministry of Science and Technology of Serbia, grant 144013 and 144018}
\address{Matemati\v{c}ki institut \\
Kneza Mihaila 36 \\
11001 Beograd, p.p. 367 \\
Serbia}
\email{\url{sasa@mi.sanu.ac.rs}}
\author[J.~D.~Phillips]{J.~D.~Phillips}
\address{Department of Mathematics \& Computer Science \\
Wabash College \\
Crawfordsville, IN 47933 USA}
\email{\url{phillipj@wabash.edu}}
\urladdr{\url{http://persweb.wabash.edu/facstaff/phillipj/}}



\begin{document}

\begin{abstract}
Right groups are direct products of  right zero semigroups and
groups and they play a significant role in the semilattice
decomposition theory of semigroups. Right groups can be
characterized as associative right quasigroups (magmas in
which left translations are bijective). If we do not assume
associativity we get right quasigroups which are not necessarily representable
as direct products of right zero semigroups and quasigroups. To
obtain such a representation, we need stronger assumptions which
lead us to the notion of \emph{right product quasigroup}. If
the quasigroup component is a (one-sided) loop, then we have a
\emph{right product (left, right) loop}.

We find a system of identities which axiomatizes right
product quasigroups, and use this to find axiom systems for
right product (left, right) loops; in fact, we can obtain each
of the latter by adjoining just one appropriate axiom to the right product
quasigroup axiom system.

We derive other properties of right product quasigroups and loops,
and conclude by showing that the axioms for right product quasigroups
are independent.
\end{abstract}

\maketitle

\section{Introduction}
\seclabel{intro}

In the semigroup literature (e.g., \cite{CliffordPreston}),
the most commonly used definition of \emph{right group}
is a semigroup $(S; \cdot)$ which is right simple
(\emph{i.e.}, has no proper right ideals) and left
cancellative (\emph{i.e.}, $xy = xz \implies y = z$).
The structure of right groups is clarified by the
following well-known representation theorem
(see \cite{CliffordPreston}):

\begin{theorem}
\thmlabel{CliffPres}
A semigroup $(S; \cdot)$ is a right group if and only if it is isomorphic to a direct
product of a group and a right zero semigroup.
\end{theorem}

There are several equivalent ways of characterizing right groups.
One of particular interest is the following: a right group is
a semigroup $(S;\cdot)$ which is also a \emph{right quasigroup},
that is, for each $a,b\in S$, there exists a unique $x\in S$ such that
$ax = b$. In a right quasigroup $(S;\cdot)$, one can define an
additional operation $\ldiv : S\times S\to S$
as follows: $z = x\ldiv y$ is the unique solution of the
equation $xz = y$. Then the following equations hold.
\[
(\eqQ{1}) \qquad x \ldiv xy = y   \qquad\qquad\qquad
(\eqQ{2}) \qquad x(x \ldiv y) = y
\]
Conversely, if we now think of $S$ as an algebra with two
binary operations then we have an equational definition.

\begin{definition}
\deflabel{rq}
An algebra $(S;\cdot, \ldiv)$ is a \emph{right quasigroup} if it satisfies
$(Q1)$ and $(Q2)$. An algebra $(S;\cdot, /)$ is a \emph{left quasigroup}
if it satisfies
\[
(\eqQ{3}) \qquad xy / y = x  \qquad\qquad\qquad
(\eqQ{4}) \qquad (x/y)y = x
\]
An algebra $(S;\cdot,\ldiv, /)$ is a \emph{quasigroup} if
it is both a right quasigroup and a left quasigroup.
\end{definition}

(We are following the usual convention that juxtaposition binds
more tightly than the division operations, which in turn bind more
tightly than an explicit use of $\cdot$. This helps avoid excessive
parentheses.)

From this point of view, a group is an associative quasigroup
with $x\ldiv y = x\inv y$ and $x\rdiv y = xy\inv$.
If $(S;\cdot,\ldiv)$ is a right group viewed as an associative
right quasigroup, then its group component has a natural right
division operation $\rdiv$. This operation can be extended to all
of $S$ as follows. We easily show that $x\ldiv x = y\ldiv y$ for all
$x,y\in S$, and then define $e = x \ldiv x$, $x^{-1} = x \ldiv e$,
and $x \rdiv y = xy^{-1}$. Note that in the right zero semigroup
component of $S$, we have $xy = x\ldiv y = x\rdiv y = y$.

If one tries to think of a right quasigroup as a ``nonassociative
right group'', one might ask if there is a representation theorem
like Theorem \thmref{CliffPres} which expresses a
right quasigroup as a direct product of a quasigroup and
a right zero semigroup. This is clearly not the case.

\begin{example}
\exmlabel{smallnotproduct}
On the set $S = \{0,1\}$, define operations
$\cdot, \ldiv : S\times S\to S$ by $x\cdot 0 = x\ldiv 0 = 1$
and $x\cdot 1 = x\ldiv 1 = 0$. Then $(S;\cdot,\ldiv)$ is a
right quasigroup which is neither
a quasigroup nor a right zero semigroup, and since
$|S| = 2$, $(S;\cdot,\ldiv)$ is also not a product of a
quasigroup and right zero semigroup.
\end{example}

For another obstruction to a representation theorem,
note that if an algebra which is a direct product
of a quasigroup and a right zero semigroup possesses a right
neutral element, then the right zero semigroup component
is trivial and the algebra is, in fact, a right loop (see below).
However, there are right quasigroups with neutral elements
which are not right loops.

\begin{example}
\exmlabel{notproduct}
Let $\bN$ be the set of natural numbers and define
\[
x \cdot y = x \ldiv y =
\begin{cases}
y & \text{if}\ x < y \\
x-y & \text{if}\ x \geq y
\end{cases}
\]
Then $(\bN; \cdot, \ldiv)$ is a right quasigroup, $0$ is a neutral element,
and $0 \cdot 1 = 1 = 2 \cdot 1$. Since $\cdot$ is not a quasigroup
operation, it follows from the preceding discussion that
$(\bN;\cdot,\ldiv)$ is not a direct product of a quasigroup and
a right zero semigroup.
\end{example}

Simply adjoining a right division operation $\rdiv$ to a right
quasigroup does not fix the problem; for instance, in either Example
\exmref{smallnotproduct} or \exmref{notproduct}, define
$x\rdiv y = 0$ for all $x,y$.

In this paper, we will investigate varieties of right quasigroups
such that there is indeed a direct product decomposition.

\begin{definition}
\deflabel{loop}
A quasigroup $(S; \cdot, \ldiv, \rdiv)$ is a $\{$\emph{left loop, right loop, loop}$\}$
if it satisfies the identity $\{\; x \rdiv x = y \rdiv y,\; x \ldiv x = y \ldiv y,\;
x \ldiv x = y \rdiv y\;\}$.

An algebra $(S; \cdot, \ldiv, \rdiv, e)$ is a \emph{pointed quasigroup}
if $(S; \cdot, \ldiv, \rdiv)$ is a quasigroup. A pointed quasigroup is a
$\{$\emph{quasigroup with an idempotent, left loop, right loop, loop}$\}$
if the distinguished element $e$ is $\{$an idempotent ($ee = e$), left neutral
($ex = x$), right neutral ($xe = x$), neutral ($ex = xe = x$)$\}$.
\end{definition}

\begin{definition}
\deflabel{language}
Let $T = \{ \cdot, \ldiv, \rdiv \}$ be the language of quasigroups and
$M$ a further (possibly empty) set of operation symbols disjoint
from $T$. The language $\hat T = T \cup M$ is \emph{an extended
language of quasigroups}.

The language $T_1 = \{ \cdot, \ldiv, \rdiv, e \}$, obtained from $T$
by the addition of a single constant, is \emph{the language of loops}.
\end{definition}

Note that we have two different algebras under the name ``loop". They are
equivalent and easily
transformed one into the other. When we need to distinguish between them
we call the algebra $(S; \cdot, \ldiv, \rdiv)$ satisfying $x \ldiv x = y \rdiv y$
``the loop in the language of quasigroups" while the algebra $(S; \cdot, \ldiv, \rdiv, e)$
satisfying identities $ex = xe = x$ is called ``the loop in the language of loops''.
Analogously we do for left and right loops.

\begin{definition} \deflabel{RPQ1}
Let $\mathbb{V}$ be a class of quasigroups.
An algebra is a \emph{right product} $\mathbb{V}$--\emph{quasigroup} if it is isomorphic
to $Q \times R$, where $Q \in \mathbb{V}$ and $R$ is a
right zero semigroup.

In particular, when $\mathbb{V}$ is the class
$\{ \Qgp, \LL, \RL, \Loop \}$ of all
$\{$\emph{quasigroups, left loops, right loops, loops}$\}$
(in the language of quasigroups)
then $\{ \RPQ, \RPLL, \RPRL, \RPL \}$ denote the class of all
right product $\mathbb{V}$--quasigroups.

If $\mathbb{V}$ is the class
$\{ \pQ, \Qi, \eQ, \Qe, \Loope \}$
of all $\{$\emph{pointed quasigroups,
quasigroups with an idempotent,
left loops, right loops, loops}$\}$ (in the language of loops), then
$\{$\RPpQ, \RPQi, \RPeQ, \RPQe, \RPLe$\}$ denote the class of all right product
$\mathbb{V}$--quasigroups.
\end{definition}

We wish to view these classes as
varieties of algebras. In order
to make sense of this, we need to adjust the type of right zero semigroups
to match that of (equational) quasigroups. We adopt the convention
suggested above.

\begin{convention}
\cnvlabel{type}
A right zero semigroup is considered to be
an algebra in $\hat T$ satisfying
$x\ldiv y = x \rdiv y = xy = y$ for all $x,y$.
\end{convention}

This convention agrees with the one used in \cite{KrapezRL, KrapezRQ}.
Different definitions of $\ldiv$ and $/$ in right zero semigroups would
affect the form of the axioms for right product quasigroups.

We also denote the class of all (pointed) right zero semigroups by $\cR \;(p\cR)$.
Then, in the language of universal algebra, the variety of all right product
 $\mathbb{V}$--quasigroups is a product $\mathbb{V} \otimes \cR$
 of independent varieties $\mathbb{V}$ and $\cR$ (see \cite{WTaylor}).

\begin{definition}
\deflabel{headtail}
If $t$ is a term, then $\{ \head(t), \tail(t)\}$ is
the $\{$first, last$\}$ variable of $t$.
\end{definition}

The following is an immediate consequence of
Definition \defref{RPQ1} and Convention \cnvref{type}.
\begin{theorem}
\thmlabel{tail}
Let $u,v$ be terms in a language extending $\{ \cdot, \ldiv, \rdiv \}$.
Then the equality $u = v$ is true in all
right product $\mathbb{V}$--quasigroups  if and only if $\tail(u) = \tail(v)$ and $u = v$ is true in all
$\mathbb{V}$--quasigroups.
\end{theorem}

In particular:

\begin{corollary}
\corlabel{cor-tail}
Let $s, t, u$ be terms in a language extending $\{ \cdot, \ldiv,
\rdiv \}$. If $s = t$ is true in all $\mathbb{V}$--quasigroups
then $s \circ u = t \circ u$ $(\circ \in \{ \cdot, \ldiv, \rdiv \})$ is
true in all right product $\mathbb{V}$--quasigroups.
\end{corollary}
We conclude this introduction with a brief discussion of
the sequel and some notation conventions.
In \S\secref{axioms}, we will consider the problem of axiomatizing
the varieties introduced by the Definition \defref{RPQ1}.
In \S\secref{prpq} we consider various properties of
right product (pointed) quasigroups and loops.
Finally, in \S\secref{indep}, we verify the independence of the axioms.

We should mention some related work by Tamura \emph{et al} and others.
\cite{Tamura61, Tamura63, Graham, Warne}.
An ``$M$-groupoid'', defined by certain axioms, turns out to be a
direct product of a right zero semigroup and a magma
with a neutral element. The axiomatic characterization of these
in \cite{Tamura61, Tamura63} is of a somewhat different character
than ours; besides the fact that they did not need to adjust
signatures since they did not consider quasigroups, their axioms
are also not entirely equational.

\begin{acknowledgment}
Our investigations were aided by the automated deduction program
\textsc{Prover9} and the finite model builder \textsc{Mace4}, both
developed by McCune \cite{McCune}
\end{acknowledgment}

\section{Axioms}
\seclabel{axioms}

We now consider the problem of axiomatizing $\RPQ$,
the class of all right product quasigroups.
One approach to axiomatization is the standard method of Knoebel
\cite{Knoebel}, which was used in \cite{KrapezRL, KrapezRQ}.
It turns out that the resulting axiom system consists of 14 identities,
most of which are far from elegant.
Another way is via independence of $\Qgp$ and $\cR$.
Using the term $\alpha(x,y) = xy \rdiv y$ (see \cite[Prop. 0.9]{WTaylor}),
we get these axioms:
\[
\begin{matrix}
xx \rdiv x = x \\
\begin{matrix}
(xy \rdiv y)(uv \rdiv v) \rdiv (uv \rdiv v) = xv \rdiv v
& (xy \cdot uv) \rdiv uv = (xu \rdiv u)(yv \rdiv v)\\
(x \ldiv y)(u \ldiv v) \rdiv (u \ldiv v) = (xu \rdiv u) \ldiv (yv \rdiv v)
& (x \rdiv y)(u \rdiv v) \rdiv (u \rdiv v) = (xu \rdiv u) \rdiv (yv \rdiv v)
\end{matrix}
\end{matrix}
\]
which we also find to be somewhat complicated.
Instead, we propose a different scheme, which
we call system (\eqA):
\begin{gather}
x \ldiv xy = y
\tag{A1}\\
x\cdot x \ldiv y = y
\tag{A2}\\
x \rdiv y\cdot y = xy \rdiv y
\tag{A3}\\
(x \rdiv y\cdot y)\rdiv z = x \rdiv z
\tag{A4} \\
xy \rdiv z\cdot z = x(y\rdiv z\cdot z)
\tag{A5}
\end{gather}

We now prove that system (\eqA) axiomatizes the variety
of right product quasigroups. It is not difficult to use the
results of \cite{KinyonPhillips} to prove this, but instead
we give a somewhat more enlightening self-contained proof.
We start with an easy observation.

\begin{lemma}
\lemlabel{quasi-semi-OK}
Every right product quasigroup satisfies system \emph{(}\eqA\emph{)}.
\end{lemma}

\begin{proof}
The quasigroup axioms (\eqQ{3}) and (\eqQ{4}) trivially
imply (\eqA{3})--(\eqA{5}), and so quasigroups satisfy (\eqA). For
each (\eqA{i}), the tails of both sides of the equation coincide.
By Theorem \thmref{tail}, we have the desired result.
\end{proof}

In an algebra $(S;\cdot,\ldiv,\rdiv)$ satisfying system (\eqA),
define a new term operation $\star : S\times S\to S$ by
\[
x \star y = xy \rdiv y = x\rdiv y\cdot y
\tag{$\star$}
\]
for all $x,y\in S$. Here the second equality follows from
(\eqA{3}), and we will use it freely without reference in what
follows.

\begin{lemma}
\lemlabel{tech}
Let $(S;\cdot,\ldiv,\rdiv)$ be an algebra satisfying system \emph{(}\eqA\emph{)}.
Then for all $x,y,z\in S$,
\begin{align}
(xy)\star z &= x(y\star z)              \eqnlabel{tech-cdot} \\
(x\ldiv y)\star z &= x \ldiv (y\star z) \eqnlabel{tech-ldiv} \\
(x\rdiv y)\star z &= x \rdiv (y\star z) \eqnlabel{tech-rdiv}
\end{align}
\end{lemma}

\begin{proof}
Equation \eqnref{tech-cdot} is just (\eqA{5}) rewritten.
Replacing $y$ with $x\ldiv y$ and using (\eqA{1}), we
get \eqnref{tech-ldiv}. Finally, for \eqnref{tech-rdiv},
we have
\begin{align*}
x \rdiv (y \star z) &= (x \star y\star z) \rdiv (y \star z)
= \lbrack (x \rdiv y\cdot y)\star z\rbrack \rdiv (y \star z) \\
&= \lbrack (x\rdiv y)(y\star z) \rbrack \rdiv (y \star z)
= (x\rdiv y)\star y\star z
= (x\rdiv y)\star z\,,
\end{align*}
using (\eqA{4}) in the first equality,
\eqnref{tech-cdot} in the third, and the rectangular
property of $\star$ in the fifth.
\end{proof}

\begin{lemma}
\lemlabel{band}
Let $(S;\cdot,\ldiv,\rdiv)$ be an algebra satisfying system \emph{(}\eqA\emph{)}.
Then $(S;\star)$ is a rectangular band.
\end{lemma}

\begin{proof}
Firstly,
\begin{equation}
\eqnlabel{band1}
(x\star y)\star z = (x \rdiv y\cdot y)\rdiv z\cdot z
= x\rdiv z\cdot z = x\star z\,,
\end{equation}
using (\eqA{4}). Replacing $x$ with $x\rdiv (y\star z)$ in
\eqnref{tech-cdot}, we get
\begin{equation}
\eqnlabel{band3}
\lbrack (x\rdiv (y\star z))y\rbrack \star z
= x\rdiv (y\star z)\cdot (y\star z)
= x\star (y\star z)\,.
\end{equation}
Thus,
\begin{equation}
\eqnlabel{band4}
x\star z = (x\star (y\star z))\star z
= (\lbrack (x\rdiv (y\star z))y\rbrack \star z)\star z
= \lbrack (x\rdiv (y\star z))y\rbrack \star z
= x\star (y\star z)\,,
\end{equation}
using \eqnref{band1}, \eqnref{band3}, \eqnref{band1} again
and \eqnref{band3} once more. Together, \eqnref{band1} and
\eqnref{band4} show that $(S;\star)$ is a semigroup
satisfying $x\star y\star z = x\star z$.

What remains is to show the idempotence of $\star$.
Replace $x$ with $x\rdiv x$ in \eqnref{tech-cdot} and set $y = z = x$, we have
\[
(x\rdiv x)(x\star x) = (x\rdiv x\cdot x)\star x
= (x\star x)\star x = x\star x\,,
\]
using \eqnref{band1}, and so
\[
x\star x = (x\rdiv x) \ldiv (x\star x)
= (x\rdiv x)\ldiv (x\rdiv x\cdot x) = x\,,
\]
using (\eqA{1}) in the first and third equalities.
\end{proof}

Let $(S;\cdot,\ldiv,\rdiv)$ be an algebra satisfying system (\eqA).
By Lemma \lemref{band}, $(S;\star)$ is a rectangular band,
and so $(S;\star)$ is isomorphic to the direct product
of a left zero semigroup and a right zero semigroup \cite{CliffordPreston}.
It will be useful to make this explicit. Introduce
translation maps in the semigroup $(S;\star)$ as follows
\[
\ell_x(y) := x \star y =: (x)r_y\,,
\]
so that the left translations $\ell_x : S\to S$ act on the left and
the right translations $r_y : S\to S$ act on the right.
Let $L = \langle \ell_x | x\in S\rangle$ and
$R = \langle r_x | x\in S\rangle$. Then $L$ is a left zero
transformation semigroup,
that is, $\ell_x \ell_y = \ell_x$, while $R$ is a right zero
transformation semigroup,
that is, $r_x r_y = r_y$. Since $\ell_x = \ell_{x\star y}$
and $r_y = r_{x\star y}$ for all $x,y\in S$, it follows easily that
the map $S\to L\times R ; x\mapsto (\ell_x,r_x)$ is an isomorphism
of semigroups.

Now we define operations $\cdot$, $\ldiv$ and $\rdiv$ on $R$
and $L$. Firstly, we define $\cdot, \ldiv, \rdiv : R\times R\to R$
by
\[
r_x \cdot r_y := r_x \ldiv r_y := r_x \rdiv r_y := r_y\,.
\]
For later reference, we formally record the obvious.

\begin{lemma}
\lemlabel{R-semi}
Let $(S;\cdot,\ldiv,\rdiv)$ be an algebra satisfying system \emph{(}\eqA\emph{)}.
With the definitions above, $(R;\cdot,\ldiv,\rdiv)$ is a right zero semigroup
\end{lemma}

It follows from Lemma \lemref{quasi-semi-OK} that $(R;\cdot,\ldiv,\rdiv)$
is an algebra satisfying system (\eqA).

\begin{lemma}
\lemlabel{R-ok}
Let $(S;\cdot,\ldiv,\rdiv)$ be an algebra satisfying system \emph{(}\eqA\emph{)}.
The mapping $S\to R ; x \mapsto r_x$ is a surjective homomorphism
of such algebras.
\end{lemma}

\begin{proof}
Firstly,
\[
(x) r_{yz} = x\star (yz) = x\star \lbrack y(z\star z)\rbrack
= x \star (yz)\star z = x\star z = xr_z
= (x)(r_y\cdot r_z)\,,
\]
using \eqnref{tech-cdot} in the third equality and $(S;\star)$
being a rectangular band in the fourth equality. Similar
arguments using \eqnref{tech-ldiv} and \eqnref{tech-rdiv}
give  $r_{y\ldiv z} = r_y \ldiv r_z$ and
$r_{y\rdiv z} = r_y \rdiv r_z$, respectively. The
surjectivity is clear.
\end{proof}

Next, we define $\cdot, \ldiv, \rdiv :
L\times L\to L$ by
\begin{align*}
(\ell_x \cdot \ell_y)(z) &= \ell_x(z) \cdot \ell_y(z) \\
(\ell_x \ldiv \ell_y)(z) &= \ell_x(z) \ldiv \ell_y(z) \\
(\ell_x \rdiv \ell_y)(z) &= \ell_x(z) \rdiv \ell_y(z)
\end{align*}
for all $x,y,z\in S$.

\begin{lemma}
\lemlabel{L-quasi}
Let $(S;\cdot,\ldiv,\rdiv)$ be an algebra satisfying system \emph{(}\eqA\emph{)}.
With the definitions above, $(L;\cdot,\ldiv,\rdiv)$ is a quasigroup.
\end{lemma}

\begin{proof}
Equations (\eqQ{1}) and (\eqQ{2}) follow immediately from
the definitions together with (\eqA{1}) and (\eqA{2}). By
(\eqA{3}), it remains to prove, say, (\eqQ{3}). For all $x,y,z\in S$,
\[
((\ell_x \cdot \ell_y) \rdiv \ell_y)(z)
= (\ell_x(z)\cdot \ell_y(z))\rdiv \ell_y(z)
= (x\star z)\star (y\star z)
= x\star z
= \ell_x(z)\,,
\]
where we have used the fact that $(S;\star)$ is a rectangular band
in the third equality.
\end{proof}

\begin{lemma}
\lemlabel{L-ok}
Let $(S;\cdot,\ldiv,\rdiv)$ be an algebra satisfying system \emph{(}\eqA\emph{)}.
The mapping $S\to L ; x \mapsto \ell_x$ is a surjective homomorphism
of such algebras.
\end{lemma}

\begin{proof}
For all $x,y,z\in S$, we compute
\begin{align*}
\ell_x(z)\cdot \ell_y(z) &= (x \star z)(y \star z)
= (x \star y \star z)(y \star z)
= \lbrack (x(y\star z)) \rdiv (y\star z)\rbrack (y \star z) \\
&= (x(y\star z)) \star y \star z
= x(y\star z \star y \star z)
= x[y \star z]
= (xy)\star z
= \ell_{xy}(z)\,,
\end{align*}
where we use rectangularity of $\star$ in the second equality,
\eqnref{tech-cdot} in the fifth, idempotence of $\star$ in the
sixth and \eqnref{tech-cdot} in the seventh. Next, if we replace
$y$ with $x\ldiv y$ and use (\eqA{1}), we get
$\ell_{x\ldiv y}(z) = \ell_x(z)\ldiv \ell_y(z)$.
Finally,
\[
\ell_x(z)\rdiv \ell_y(z) = (x\star z)\rdiv (y\star z)
= ((x\star z)\rdiv y)\star z = (x\rdiv y)\star z
= \ell_{x\rdiv y}(z)\,,
\]
using \eqnref{tech-rdiv} in the second equality and
(\eqA{5}) in the third.
\end{proof}

We now turn to the main result of this section.

\begin{theorem}
\thmlabel{repr}
An algebra $(S;\cdot,\ldiv,\rdiv)$ is a right product quasigroup
if and only if it satisfies \emph{(}\eqA\emph{)}.
\end{theorem}

\begin{proof}
The necessity is shown by Lemma \lemref{quasi-semi-OK}.
Conversely, if $(S;\cdot,\ldiv,\rdiv)$ satisfies (\eqA), then by
Lemmas \lemref{R-ok} and \lemref{L-ok}, the mapping
$S\to L\times R ; x \mapsto (\ell_x, r_x)$ is a surjective
homomorphism. This map is, in fact, bijective, since as already
noted, it is an
isomorphism of rectangular bands. By Lemmas \lemref{R-semi}
and \lemref{L-quasi}, $L\times R$ is a right product quasigroup,
and thus so is $S$.
\end{proof}

\begin{remark}
There are other choices of axioms for right product quasigroups.
For instance, another system equivalent to (\eqA) consists of
(\eqA{1}), (\eqA{2}), (\eqA{3}) and the equations
\[
(\eqB{1}) \quad xx \rdiv x = x   \qquad\qquad\qquad
(\eqB{2}) \quad (xy\cdot (z\rdiv u)) \rdiv (z\rdiv u) = x(yu\rdiv u)\,.
\]
Call this system (\eqB). We omit the proof of the equivalence of
systems (\eqA) and (\eqB). One can use the results of
\cite{KinyonPhillips} to prove the system (\eqB) variant of
Theorem \thmref{repr} as follows: (\eqA{1}) and (\eqA{2})
trivially imply the equations
\[
(\eqA{3}^{\prime})\quad
x (x \ldiv y) = x \ldiv xy \qquad
(\eqB{1}^{\prime})\quad
x \ldiv xx = x \qquad
(\eqB{2}^{\prime})\quad (x\ldiv y) \ldiv ((x\ldiv y)\cdot zu) = (x \ldiv xz) u\,.
\]
By \cite{KinyonPhillips}, (\eqA{3}), (\eqA{3}$^{\prime}$),
(\eqB{1}), (\eqB{1}$^{\prime}$), (\eqB{2}$^{\prime}$) and
(\eqB{2}$^{\prime}$) axiomatize
the variety of \emph{rectangular quasigroups}, each of which is
a direct product of a left zero semigroup, a quasigroup and a
right zero semigroup. By (\eqA{1}) and (\eqA{2}), the left zero
semigroup factor must be trivial, and so a system satisfying
system (\eqB) must be a right product quasigroup.
\end{remark}

We conclude this section by considering other
varieties of right product quasigroups.
Utilizing \cite{Krapez:V} we get:

\begin{theorem}
\thmlabel{RPV}
Let $\mathbb{V}$ be a variety of quasigroups axiomatized by additional
identities:
\[
s_i = t_i \tag{\eqV$_i$}
\]
 $(i \in I)$ in an extended language $\hat T$ and let $z$ be a
 variable which does not occur in any $s_i, t_i$.
Then the variety \RPV\; of right product $\mathbb{V}$--quasigroups
can be axiomatized by system (\eqA)
together with (for all $i \in I$):
\begin{equation}
s_i z = t_i z \tag{{\eqhV}$_i$}\,.
\end{equation}
\end{theorem}

\begin{proof}
Both $\mathbb{V}$--quasigroups and right zero semigroups satisfy
system (\eqA) and all ({\eqhV}$_i$), $i \in I$, and thus
so do their direct products i.e. right product
$\mathbb{V}$--quasigroups.

Conversely, if an algebra satisfies system (\eqA), it is
 a right product quasigroup by Theorem \thmref{repr}. Since
 all ({\eqhV}$_i$) are satisfied, the quasigroup factor
 has to satisfy them, too.
 But in quasigroups, the identities ({\eqhV}$_i$)
 are equivalent to the identities ({\eqV}$_i$) and these
 define the variety $\mathbb{V}$.
\end{proof}

\begin{theorem}
\thmlabel{RPV2}
Theorem \thmref{RPV} remains valid if we replace \emph{(}{\eqhV}$_i$\emph{)}
by any of the following families of identities:
%(provided $z$ does not occur in $s_i, t_i$):
\begin{gather*}
s_i \ldiv z = t_i \ldiv z \\
s_i \rdiv z = t_i \rdiv z \\
z \rdiv (s_i \ldiv z) = (z \rdiv t_i) \ldiv z \\
s_i = (t_i \cdot \tail (s_i)) \rdiv \tail (s_i) \\
s_i = (t_i \rdiv \tail (s_i)) \cdot \tail (s_i) \\
s_i = t_i \;\;\; (\text{if } \tail (s_i) = \tail (t_i)).
\end{gather*}
\end{theorem}

\begin{example}
\exmlabel{right groups}
Adding associativity $x \cdot yz = xy \cdot z$ to system
(\eqA) gives yet another axiomatization of
\emph{right groups}.
\end{example}

\begin{example}
\exmlabel{RP comm. Q}
Right product commutative quasigroups are right product quasigroups
satisfying $xy \cdot z = yx \cdot z$.
However,  commutative right product quasigroups are just
commutative quasigroups.
\end{example}

Obviously:

\begin{corollary}
\corlabel{V to RPV}
If the variety $\mathbb{V}$ of quasigroups is defined by the
identities $s_i = t_i$ $((i \in I)$
such that $\tail(s_i) = \tail(t_i)$ for all $i \in I$,
then the class of all right product quasigroups satisfying
identities $s_i = t_i (i \in I)$ is the class of all
right product $\mathbb{V}$--quasigroups.

If $\tail(s_i) \neq \tail(t_i)$ for some $i \in I$,
then the class of all right product quasigroups satisfying
identities $s_i = t_i$ $(i \in I)$ is just the class of all
$\mathbb{V}$--quasigroups.
\end{corollary}

\begin{example}
\exmlabel{RPpQ}
The variety \RPpQ\; is defined by adding a constant to the language
of quasigroups, not by any extra axioms.
\end{example}

\begin{example}
\exmlabel{RPQi}
The variety \RPQi\; of all \emph{right product quasigroups with an idempotent}
may be axiomatized by system (\eqA) and $ee = e$.
\end{example}

\begin{corollary}
\corlabel{RPLL}
A right product quasigroup is a right product left loop
iff it satisfies any (and hence all) of the following axioms:
\begin{gather}
(x \rdiv x)y = y \tag{LL1} \eqnlabel{LL} \\
(x \rdiv x)z = (y \rdiv y)z \eqnlabel{LL2} \tag{LL2} \\
(x \circ y) \rdiv (x \circ y) = y \rdiv y \eqnlabel{LL6} \tag{LL3}
\end{gather}
where $\circ$ is any of the operations $\cdot, \ldiv, \rdiv$.
\end{corollary}

\begin{proof}
In a quasigroup, identities \eqnref{LL}, \eqnref{LL2} and
\eqnref{LL6} are equivalent to each other and to $x\rdiv x = y\rdiv y$,
and so a quasigroup satisfying either axiom is a left loop.
Conversely, in a left loop with left neutral element $e$, we have
$e = x\rdiv x$, and so \eqnref{LL}, \eqnref{LL2} and \eqnref{LL6}
hold. Thus a quasigroup satisfies either (and hence all) of
\eqnref{LL}, \eqnref{LL2}, \eqnref{LL6} if and only if it is a left loop.

On the other hand, \eqnref{LL}, \eqnref{LL2} and \eqnref{LL6}
trivially hold in right zero semigroups by Convention
\cnvref{type}. Putting this together, we have the desired result.
\end{proof}

In the language of loops we have:

\begin{corollary}
\corlabel{RPeQ}
A right product quasigroup
is a right product left loop if and only if it satisfies
the identity $ex = x$.
\end{corollary}

Similarly:

\begin{corollary}
\corlabel{RPRL}
A right
product quasigroup is a right product right loop if and only if it satisfies
any (and hence all) of the following axioms:
\begin{gather*}
x(y \ldiv y) \cdot z = xz \\
(x \ldiv x)z = (y \ldiv y)z \\
(x \circ y) \ldiv (x \circ y) = y \ldiv y
\end{gather*}
where $\circ$ is any of the operations $\cdot, \ldiv, \rdiv$.
\end{corollary}

\begin{corollary}
\corlabel{RPQe}
A right product quasigroup is a right product right loop
(in the language of loops)
iff it satisfies the identity $xe \cdot y = xy$.
\end{corollary}

\begin{corollary}
\corlabel{RPL}
A right
product quasigroup is a right product loop if and only if it satisfies
any (and hence all) of the following axioms:
\begin{gather}
(x \ldiv x)y = y \tag{L} \\
x(y \rdiv y) = xy \rdiv y \notag \\
x(y \rdiv y) = (x \rdiv y)y \notag \\
(x \ldiv x)z = (y \rdiv y)z \notag \\
(x \circ y) \ldiv (x \circ y) = y \rdiv y \notag \\
(x \circ y) \rdiv (x \circ y) = y \ldiv y \notag
\end{gather}
where $\circ$ is any of the operations $\cdot, \ldiv, \rdiv$.
\end{corollary}

\begin{corollary}
\corlabel{RPQ1}
A right product quasigroup is a right product loop
(in the language of loops)
iff it satisfies both $ex = x$ and $xe \cdot y = xy$.
\end{corollary}

\section{Properties of Right Product (Pointed) Quasigroups}
\seclabel{prpq}

Calling upon the tools of universal algebra, we now examine
some properties of right product quasigroups.
We will use the following standard notation.

\begin{definition}
\deflabel{standard}
\hfill

\noindent
\begin{tabular}{lcl}
$E_S$ & -- & the subset of all idempotents of $S$. \\
$\Sub(S)$ & -- & the lattice of all subalgebras of $S$. \\
$\Sub^0(S)$ & -- & the lattice of all subalgebras of $S$ with the empty set adjoined as the \\
&& smallest element (used when two subalgebras have an empty intersection). \\
$\Con(S)$ & -- & the lattice of all congruences of $S$. \\
$\Eq(S)$ & -- & the lattice of all equivalences of $S$. \\
$\Hom(S,T)$ & -- & the set of all homomorphisms from $S$ to $T$. \\
$\End(S)$ & -- & the monoid of all endomorphisms of $S$. \\
$\Aut(S)$ & -- & the group of all automorphisms of $S$. \\
$\Free({\mathbb{V}}, n)$ & -- & the free algebra with $n$ generators in the variety $\mathbb{V}$. \\
$\Var({\mathbb{V}})$ & -- & the lattice of all varieties of a class ${\mathbb{V}}$ of algebras.
\end{tabular}
\smallskip

\noindent In addition, $R_n$ will denote the unique $n$-element
right zero semigroup -- which also happens to be free.
However, note that in the language of loops, the free right zero semigroup
generated by $n$ elements is $R_{n+1}$.
\end{definition}

\subsection{The word problem}
\seclabel{word}
\hfill{ }
\smallskip

Using a well-known result of Evans \cite{Evans} we have the
following corollary of Theorem \thmref{tail}:

\begin{corollary}
\corlabel{word}
The word problem for right product $\mathbb{V}$--quasigroups
is solvable if and only if it is solvable for $\mathbb{V}$--quasigroups.
\end{corollary}

In particular:
\begin{corollary}
\corlabel{word:RPQ}
The word problem for $\{\RPQ, \RPLL, \RPRL, \RPL \}$\; is solvable.
\end{corollary}

Likewise:

\begin{corollary}
\corlabel{word:RPpQ}
The word problem for $\{\RPpQ, \RPQi, \RPeQ, \RPQe, \RPLe \}$\; is solvable.
\end{corollary}

\subsection{Properties of right product quasigroups and loops}
\seclabel{properties}
\hfill{ }
\smallskip

The following corollaries are special cases of results in universal algebra
(see \cite{WTaylor}).

\begin{corollary}
\corlabel{RPQfunctors}
For all $Q, Q' \in \Qgp$ and $R, N \in \cR$:
\begin{enumerate}[label=\emph{(}\arabic*\emph{)}]
\item $E_{Q \times R} = E_Q \times R$ \\
      in particular:
      \begin{itemize}
      \item[-] $Q \times R$ have idempotents if and only if $Q$ have them
      \item[-] $E_{Q \times R}$ is subalgebra of $Q \times R$ if and only if
      $E_Q$ is subalgebra of $Q$
      \item[-] $Q \times R$ is a groupoid of idempotents if and only if $Q$ is
      \end{itemize}
\item $\Sub^0(Q \times R)= (\Sub(Q)\times
    ({\bf 2}^R\smallsetminus\{\emptyset\}))\cup\{\emptyset\}$
\item $\Con(Q \times R) = \Con(Q) \times \Eq(R)$
\item $\Hom(Q \times R,\,Q' \times N) = \Hom(Q,Q') \times N^R$
\item $\End(Q \times R) = \End(Q) \times R^R$
\item $\Aut(Q \times R) = \Aut(Q) \times S_{|R|} .$
\end{enumerate}
\end{corollary}

Having a distinguished element changes the properties of a variety radically.
For example,
if $e = (i,j)$ is a distinguished element of the right product pointed quasigroup
$Q \times R$ then there is always the smallest subalgebra $<i> \times \{j\}$.
So, in case of right product pointed quasigroups, the results analogous to
Corollary \corref{RPQfunctors} are actually somewhat different in character.

\begin{corollary}
\corlabel{RPpQfunctors}
For all $Q, Q' \in \pQ$\; and $R, N \in \cR$\; with a distinguished element $j$:
\begin{enumerate}[label=\emph{(}\arabic*\emph{)}]
\item $E_{Q \times R} = E_Q \times R$ \\
      in particular:
      \begin{itemize}
      \item[-] $Q \times R$ has idempotents if and only if $Q$ has them.
      \item[-] $E_{Q \times R}$ is subalgebra of $Q \times R$ if and only if $E_Q$ is subalgebra of $Q$.
      \item[-] $Q \times R$ is a groupoid of idempotents if and only if $Q$ is.
      \end{itemize}
\item $\Sub(Q \times R)= \Sub(Q)\times\{ Y \subseteq R \mid j \in Y \}
                       \simeq \Sub(Q) \times {\bf 2}^{R\smallsetminus\{j\}}$
\item $\Con(Q \times R) = \Con(Q) \times \Eq(R)$
\item $\Hom(Q \times R,\,Q' \times N) = \Hom(Q,Q') \times \{f : R \rightarrow N \mid f(j)=j\}
                                      \simeq \Hom(Q,Q') \times N^{R\smallsetminus \{j\}}$
\item $\End(Q \times R) = \End(Q) \times \{ f : R \rightarrow R \mid f(j)=j \}
                        \simeq \End(Q) \times R^{R \smallsetminus \{j\}}$
\item $\Aut(Q \times R) = \Aut(Q) \times \{ f \in S_R \mid f(j) = j \}
                        \simeq \Aut(Q) \times S_{|R|-1} .$
\end{enumerate}
\end{corollary}

\begin{corollary}
\corlabel{free}
If $\mathbb{V}$ is one of $\Qgp, \LL, \RL, \Loop$, then
$\Free(\RPV,n)\simeq \Free(\mathbb{V}, n)\times R_n$.
\end{corollary}

\begin{corollary}
\corlabel{pfree}
If $\mathbb{V}$ is one of $\pQ, \Qi, \eQ, \Qe, \Loope$, then
$\Free(\RPV, n)\simeq \Free(\mathbb{V}, n)\times R_{n+1}$.
\end{corollary}

\begin{corollary}
\corlabel{var}
If $\mathbb{V}$ is one of the above varieties of (pointed) quasigroups, then
$\Var(\RPV) \simeq \Var({\mathbb{V}}) \times {\bf 2}$.
\end{corollary}

All cases suggested by Corollary \corref{RPQfunctors}(1) can actually
occur. In the examples below,
right product quasigroups are in fact quasigroups and thus we display
the Cayley tables of the multiplication only.

\begin{example}
\exmlabel{noid}
Tables \tabref{noid} give a right product quasigroup with no idempotents
(on the left) and an idempotent right product quasigroup (on the right).
\begin{table}[htb]
\[
\begin{array}{c|ccc}
\cdot &  0& 1& 2\\
\hline
0 &  1 &  0& 2\\
1 &  0 &  2& 1\\
2 &  2 &  1& 0\\
\end{array}
\qquad
\begin{array}{c|ccc}
\cdot &  0& 1& 2\\
\hline
0 &  0 &  2& 1\\
1 &  2 &  1& 0\\
2 &  1 &  0& 2\\
\end{array}
\]
\caption{A right product quasigroup with no idempotents and and idempotent
right product quasigroup}
\end{table}
\tablabel{noid}
\end{example}

\begin{example}
\exmlabel{not closed}
Tables \tabref{notclosed} give a right product quasigroup in which
$E_S$ is not a subalgebra (on the left) and a right product quasigroup
in which $E_S$ is a nontrivial subalgebra (on the right).
\begin{table}[htb]
\[
\begin{array}{lcr}
% multiplication
\begin{array}{c|cccc}
\bullet &  0& 1& 2& 3\\
\hline
0 &  0 &  2& 1 & 3\\
1 &  3 &  1& 2 & 0\\
2 &  1 &  3& 0 & 2\\
3 &  2 &  0& 3 & 1\\
\end{array}
&& \qquad
\begin{array}{c|cccccc}
\bullet &  0& 1& 2 & 3 & 4 & 5\\
\hline
0 &  0 &  2& 1 & 3 & 5 & 4\\
1 &  2 &  1& 0 & 5 & 4 & 3\\
2 &  1 &  0& 2 & 4 & 3 & 5\\
3 &  3 &  5& 4 & 0 & 2 & 1\\
4 &  5 &  4& 3 & 2 & 1 & 0\\
5 &  4 &  3& 5 & 1 & 0 & 2\\
\end{array}
\end{array}
\]
\caption{$E_S$ is not closed, $E_S$ is a nontrivial subalgebra}
\tablabel{notclosed}
\end{table}
\end{example}

Moreover, we have the following. These are immediate consequences
of well understood properties of quasigroups and semigroups.

\begin{theorem}
\thmlabel{RPQ:max}
Let $S = Q \times R$ be a right product quasigroup. Then:
\begin{enumerate}[label=\emph{(}\arabic*\emph{)}]
\item If $Q_m \enskip (m \in M)$ is the (possibly empty) family
      of all maximal subquasigroups of $Q$ then
      $Q_m \times R \enskip (m \in M), \enskip
      Q \times (R \setminus \{ r \}) \enskip (r \in R)$
      is the family of all maximal right product subquasigroups of $S$.
\item There are $|R|$ maximal subquasigroups of $S$.
      They are all mutually isomorphic and of the form
      $Q \times \{ r \} \; (r \in R)$.
\\ From now on we assume $E_S \neq \emptyset \enskip (\text{i.e. } E_Q \neq \emptyset)$.
\item If $E_S$ is subalgebra then it is the largest subalgebra of
      idempotents of $S$.
\item There are $|E_S|$ maximal left zero subsemigroups of $S$.
      They are all singletons $\{ e \} \; (e \in E_S)$.
\item There are $|E_Q|$ maximal right zero subsemigroups of $S$.
      They are all mutually isomorphic and of the form
      $\{ i \} \times R \; (i \in E_Q)$.
\\ From now on we assume $E_Q = \{ i \}$.
\item $E_S = \{i\} \times R$ is the unique largest subband of $S$,
      which happens to be a right zero semigroup.
\item If the quasigroup $Q$ is a left loop then the left neutral $i$ of $Q$
      is the only idempotent of $Q$ and:
      \begin{itemize}
      \item[-] $E_S = \{ a \rdiv a \; | \; a \in S \}$.
      \item[-] The element $e \in S$ is a left neutral if and only if it is an idempotent.
      \item[-] The maximal subquasigroups
               $Q \times \{r\} = S e = \{ x \in S \; | \; x \rdiv x = e \}
               \enskip (e \in E_S, e = (i,r))$
               are maximal left subloops of $S$.
      \item[-] For all $e \in E_S \enskip S \simeq Se \times E_S$
               and the isomorphism is given by $f(x) = (xe \rdiv e, x \rdiv x)$.
      \end{itemize}
\item If the quasigroup $Q$ is a right loop then the right neutral $i$ of $Q$
      is the only idempotent of $Q$ and:
      \begin{itemize}
      \item[-] $E_S = \{ a \ldiv a \; | \; a \in S \}$.
      \item[-] $S$ has a right neutral if and only if $|R| = 1$ and then the right
               neutral is unique. In this case $S$ is a right loop.
      \item[-] The maximal subquasigroups
               $Q \times \{r\} = S e = \{ x \in S \; | \; x \ldiv x = e \}
               \enskip (e \in E_S, e = (i,r))$
               are maximal right subloops of $S$.
      \item[-] For all $e \in E_S \enskip S \simeq Se \times E_S$
               and the isomorphism is given by $f(x) = (xe, x \ldiv x)$.
      \end{itemize}
 \item If the quasigroup $Q$ is a loop then: %the neutral element $i$ of $Q$
%      is the only idempotent of $Q$ and:
      \begin{itemize}
      \item[-] The element $e \in S$ is a left neutral if and only if it is an idempotent.
      \item[-] $S$ has a right neutral if and only if $|R| = 1$ and then the right
               neutral is unique. In this case $S$ is a loop.
      \item[-] The maximal subquasigroups
               $Q \times \{r\} = S e = \{ x \in S \; | \; x \ldiv x = x \rdiv x = e \}
               \enskip (e \in E_S, e = (i,r))$
               are maximal subloops of $S$.
      \item[-] For all $e \in E_S \enskip S \simeq Se \times E_S$
               and the isomorphism is given by $f(x) = (xe, x \rdiv x)$.
      \end{itemize}
\end{enumerate}
\end{theorem}

\begin{theorem}
\thmlabel{RPpQ:max}
Let $S = Q \times R$ be a right product pointed quasigroup
with a distinguished element $e = (i,j)$. Then:
\begin{enumerate}[label=\emph{(}\arabic*\emph{)}]
\item If $Q_m \enskip (m \in M)$ is the (possibly empty) family
      of all maximal pointed subquasigroups of $Q$ then
      $Q_m \times R \enskip (m \in M), \enskip
      Q \times (R \setminus \{ r \}) \enskip (r \in R \setminus \{ j \})$
      is the family of all maximal right product pointed subquasigroups of $S$.
\item $S e = Q \times \{ j \}$ is the largest pointed subquasigroup of $S$.
\\ From now on we assume $E_S \neq \emptyset \enskip (\text{i.e. } E_Q \neq \emptyset)$.
\item If $E_S$ is subalgebra then it is the largest subalgebra of
      idempotents of $S$.
\item $S e \cap E_S = \{ e \}$ is the largest pointed left zero subsemigroup of $S$
      if and only if $e \in E_S$.
\item $E_S = \{ i \} \times R$ is the largest pointed right zero subsemigroup of $S$
      if and only if $i \in E_Q$.
\item $S \simeq Se \times E_S$ and the isomorphism is given by
      $f(x) = (xe \rdiv e, ex \rdiv x)$.
\\ From now on we assume $E_Q = \{ i \}$.
\item $E_S = \{i\} \times R$ is the unique largest pointed subband of $S$,
      which happens to be a pointed right zero semigroup.
\item If the element $i$ is the left neutral of $Q$, it
      is the only idempotent of $Q$ and:
      \begin{itemize}
      \item[-] $E_S = \{ a \rdiv a \; | \; a \in S \}$.
      \item[-] The element $a \in S$ is a left neutral if and only if it is an idempotent.
      \item[-] $S e = \{ x \in S \; | \; x \rdiv x = e \}$
               is the largest left subloop of $S$.
      \item[-] The isomorphism $S \simeq Se \times E_S$
               is given by $f(x) = (xe \rdiv e, x \rdiv x)$.
      \end{itemize}
\item If the element $i$ is the right neutral of $Q$, it
      is the only idempotent of $Q$ and:
      \begin{itemize}
      \item[-] $E_S = \{ a \ldiv a \; | \; a \in S \}$.
      \item[-] $S$ has a right neutral if and only if $|R| = 1$ and then the right
               neutral is $e$. In this case $S$ is a right loop.
      \item[-] $S e = \{ x \in S \; | \; x \ldiv x = e \}$
               is the largest right subloop of $S$.
      \item[-] The isomorphism $S \simeq Se \times E_S$
               is given by $f(x) = (xe, x \ldiv x)$.
      \end{itemize}
 \item If the element $i$ is the two--sided neutral of $Q$, it
      is the only idempotent of $Q$ and:
      \begin{itemize}
      \item[-] The element $a \in S$ is a left neutral if and only if it is an idempotent.
      \item[-] $S$ has a right neutral if and only if $|R| = 1$ and then the right
               neutral is $e$. In this case $S$ is a loop.
      \item[-] $S e = \{ x \in S \; | \; x \ldiv x = x \rdiv x = e \}$
               is the largest subloop of $S$.
      \item[-] The isomorphism $S \simeq Se \times E_S$
               is given by $f(x) = (xe, x \rdiv x)$.
      \end{itemize}
\end{enumerate}
\end{theorem}

\subsection{The equation xa=b}
\label{xa = b}
\hfill{ }
\smallskip

Since a right product quasigroup is a right quasigroup the equation
$ax = b$ has the unique solution $x = a \ldiv b$. For the equation
$xa = b$, the situation is not so clear.

We solve the equation $xa = b$ using the notion of \emph{reproductivity}.
The related notion of \emph{reproductive general solution} was
defined by E. Schr\"{o}der \cite{Schroder} for Boolean equations
and studied by L. L\"{o}wenheim \cite{Lowenheim08, Lowenheim10}
who also introduced
the term ``reproductive". More recently, S. B. Pre\v{s}i\'c
made significant contributions to the notion of reproductivity
\cite{Presic68,Presic72, Presic00}. For an introduction to reproductivity,
see S. Rudeanu \cite{Rudeanu}.

\begin{definition}
\deflabel{reproductive}
Let $S \neq \emptyset$ and $F : S \longrightarrow S$. The equation
$x = F(x)$ is \emph{reproductive} if for all $x \in S \;\; F(F(x)) = F(x)$.
\end{definition}
The most significant properties of reproductivity are:
\begin{theorem}
\thmlabel{general solution}
A general solution of the reproductive equation $x = F(x)$ is given by:
$x = F(p) \;\; (p \in S)$.
\end{theorem}
\begin{theorem}[S. B. Pre\v{s}i\'c]
\thmlabel{Presic}
Every consistent equation has an equivalent reproductive equation.
\end{theorem}

We now apply these results to our equation $xa = b$.

\begin{theorem}
\thmlabel{xa=b}
\begin{enumerate}[label=\emph{(}\arabic*\emph{)}]
\item In right product
quasigroups, the equation $xa = b$ is consistent  if and only if $(b / a)a = b$.
\item In right product $\{$ left, right $\}$ loops, the consistency of $xa = b$
is equivalent to $\{ a \rdiv a = b \rdiv b, a \ldiv a = b \ldiv b \}$.
\item If the equation $xa = b$ consistent, then it is equivalent to the reproductive
equation $x = (b \rdiv a)x \rdiv x$, and thus its general solution is given by
$x = (b / a)p / p$  ($p \in S$). There are exactly $|R|$ distinct solutions.
\item If a right product quasigroup $S = Q \times R$ has idempotents, then the general
solution of the consistent equation $xa = b$ may be given in the form
$x = (b \rdiv a)e \rdiv e \quad (e \in E_S)$.
\item If the quasigroup $Q$ has a unique idempotent, then any idempotent $e \in E_S$
defines the unique solution $x = (b \rdiv a)e \rdiv e$ of $xa = b$.
\item In a right product right loop, the
general solution of $xa = b$ may be simplified to
$x = (b \rdiv a) e \quad (e \in E_S)$.
\end{enumerate}
\end{theorem}

\begin{proof}
(1) If the equation is consistent then there is at least one solution
$x = c$. It follows that $b = ca$ and
$(b \rdiv a)a = (ca \rdiv a)a = ca = b$. If we assume that
$(b \rdiv a)a = b$ then $x = b / a$ is one solution of the equation
$xa = b$, which therefore must be consistent.

(2) Assume $S$ is a right product left loop.
If $xa = b$ is consistent then $b \rdiv b = xa \rdiv xa = a \rdiv a$.
For the converse assume $S = Q \times R$ for some left loop $Q$ with
the left neutral $i$ and a right zero semigroup $R$. Let $a = (a_1, a_2)$
and $b = (b_1, b_2)$.
Then $(i, a_2) = (a_1 \rdiv a_1, a_2 \rdiv a_2) = a \rdiv a = b \rdiv b =
(b_1 \rdiv b_1, b_2 \rdiv b_2) = (i, b_2)$ i.e. $a_2 = b_2$
which is equivalent to the consistency of $xa = b$.

Now assume $S$ is a right product right loop.
If $xa = b$ is consistent then $b \ldiv b = xa \ldiv xa = a \ldiv a$.
For the converse assume $S = Q \times R$ for some right loop $Q$
with the right neutral $i$ and a right zero semigroup $R$. Let $a
= (a_1, a_2)$ and $b = (b_1, b_2)$.
Then $(i,a_2) = (a_1 \ldiv a_1, a_2 \ldiv a_2) = a \ldiv a = b \ldiv b
= (b_1 \ldiv b_1, b_2 \ldiv b_2) = (i, b_2)$ i.e. $a_2 = b_2$
which is equivalent to the consistency of $xa = b$.

(3) Let the equation $xa = b$ be consistent.
Then $(b \rdiv a)x \rdiv x = (xa \rdiv a)x \rdiv x = xx \rdiv x = x$.
Conversely, if $x = (b \rdiv a)x \rdiv x$ then
$xa = ((b \rdiv a)x \rdiv x)a = (b \rdiv a)a = b$. Therefore,
equations $xa = b$ and $x = (b \rdiv a)x \rdiv x$ are equivalent.
The form of the later equation is $x = F(x)$ where
$F(x) = (b \rdiv a)x \rdiv x$. Also,
$F(F(x)) = ((b \rdiv a) \cdot F(x)) \rdiv F(x) =
((b \rdiv a) ((b \rdiv a)x \rdiv x)) \rdiv ((b \rdiv a)x \rdiv x) =
(b \rdiv a)x \rdiv x = F(x)$ so equation $x = F(x)$ is reproductive.
Its general solution is $x = F(p) = (b \rdiv a)p \rdiv p  (p \in S)$.

Without loss of generality we may assume that $S$ is $Q \times R$
for some quasigroup $Q$ and a right zero semigroup $R$. Let
$a = (a_1, a_2), b = (b_1, b_2), x = (x_1, x_2)$ and
$p = (p_1, p_2)$.
The consistency of $xa = b$ reduces to $(b_1, b_2) = (b \rdiv a)a =
((b_1 \rdiv a_1)a_1, a_2) = (b_1, a_2)$ i.e. $a_2 = b_2$.
In that case,
the solutions of $xa = b$ are $x = (b \rdiv a)p \rdiv p =
((b_1 \rdiv a_1)p_1 \rdiv p_1, (b_2 \rdiv a_2)p_2 \rdiv p_2) =
(b_1 \rdiv a_1, p_2)$. Evidently,
the number of different solutions of $xa = b$  is $|R|$.

(4) Let $S = Q \times R$ and let $i$ be an idempotent of $Q$.
For every $p = (p_1, p_2)$ there is an idempotent $e = (i, p_2)$
of $S$ such that $x = (b \rdiv a)p \rdiv p = (b_1 \rdiv a_1, p_2) =
((b_1 \rdiv a_1)i \rdiv i, (b_2 \rdiv a_2)p_2 \rdiv p_2) = (b \rdiv a)e \rdiv e$.

(5) If the idempotent $i \in Q$ is unique then $E_S$ has exactly
$|E_S| = |R|$ idempotents, just as many as the equation $xa = b$
has solutions.

(6) Let $S$ be a right product right loop and $e = (i,r)$.
Then
\begin{align*}
x &= (b \rdiv a)e \rdiv e = ((b_1 \rdiv a_1)i \rdiv i, (b_2 \rdiv a_2)r \rdiv r)
= (b_1 \rdiv a_1, r) \\
&= ((b_1 \rdiv a_1)i, (b_2 \rdiv a_2)r) =
(b_1 \rdiv a_1, b_2 \rdiv a_2)(i, r) = (b \rdiv a)e\,.
\end{align*}
\end{proof}

\subsection{Products of sequences of elements including idempotents}
\seclabel{products}
\hfill{ }
\smallskip

We use $\varrho (a_i, a_{i+1},\ldots ,a_j)$ to denote
\emph{the right product} i.e. the product of $a_i, \ldots ,a_j$
with brackets associated to the right. More formally,
$\varrho(a_i) = a_i \quad (1 \leq i \leq n)$ and
$\varrho(a_i, a_{i+1},\ldots ,a_j) = a_i \cdot \varrho(a_{i+1},\ldots ,a_j)
\quad (1 \leq i < j \leq n)$.

Further, we define $\varrho(\pm a_n) = a_n$ and
\[
\varrho(a_i, a_{i+1},\ldots ,a_j, \pm a_n) =
    \begin{cases}
    \varrho(a_i,\ldots ,a_j); & \text{if $j = n$}  \\
    \varrho(a_i,\ldots ,a_j, a_n); & \text{if $j < n$.}
    \end{cases}
\]
In short, $a_n$ should appear in the product $\varrho(a_i, \ldots, a_j, \pm a_n)$,
but only once.

The following is an analogue of Theorem 2.4 from \cite{KrapezRL}.

\begin{theorem}
\thmlabel{rightproduct}
Let $a_1,\ldots,a_n \quad (n > 0)$
be a sequence of elements of the right product left loop
$S$, such that $a_{p_1}, \ldots , a_{p_m} \quad (1 \leq p_1 < \ldots < p_m \leq n;
\quad 0 \leq m \leq n)$ comprise exactly the idempotents among $a_1, \ldots, a_n$.
Then $\varrho(a_1,\ldots,a_n) = \varrho(a_{p_1},\ldots,a_{p_m},\pm a_n)$.
\end{theorem}

\begin{proof}
The proof is by induction on $n$.

\noindent (1) $n = 1$.

\noindent If $m = 0$ then $\varrho(a_1) = a_1 = \varrho(\pm a_1)$.

\noindent If $m = 1$ then $\varrho(a_1) = a_1 = \varrho(a_1, \pm a_1)$.

\noindent (2) $n > 1$.

\noindent If $a_1 \in E_S$ (i.e. $1 < p_1$) then
$\varrho(a_1,\ldots,a_n) = a_1 \cdot \varrho(a_2,\ldots,a_n) = \varrho(a_2,\ldots,a_n)$
which is equal to $\varrho(a_{p_1},\ldots,a_{p_m},\pm a_n)$
by the induction argument.

\noindent If $a_1 \notin E_S$ (i.e. $1 = p_1$) then, using induction argument again,
we get
$\varrho(a_1,\ldots,a_n) = a_1 \cdot \varrho(a_2,\ldots,a_n) =
a_{p_1} \cdot \varrho(a_{p_2},\ldots,a_{p_m}, \pm a_n) =
\varrho(a_{p_1},a_{p_2},\ldots,a_{p_m}, \pm a_n)$.
\end{proof}

Analogously to $\varrho(\ldots)$,
we use $\lambda(a_i, \ldots, a_{j-1}, a_j)$ to denote
\emph{the left product} i.e. the product of $a_i, \ldots ,a_j$
with brackets associated to the left. Formally,
$\lambda(a_i) = a_i \quad (1 \leq i \leq n)$ and
$\lambda(a_i, \ldots, a_{j-1}, a_j) = \lambda(a_i, \ldots ,a_{j-1}) \cdot a_j
\quad (1 \leq i < j \leq n)$.

Further, we define $\lambda(\pm a_1, \pm a_n) = a_1$ if $n = 1$,
$\lambda(\pm a_1, \pm a_n) = a_1 a_n$ if $n > 1$ and
\[
\lambda(\pm a_1, a_i, \ldots ,a_j, \pm a_n) =
    \begin{cases}
    \lambda(a_i,\ldots ,a_j); & \text{if $1 = i \leq j = n$} \\
    \lambda(a_i,\ldots ,a_j, a_n); & \text{if $1 = i \leq j < n$} \\
    \lambda(a_1, a_i, \ldots, a_j); & \text{if $1 < i \leq j = n$} \\
    \lambda(a_1, a_i, \ldots, a_j, a_n); & \text{if $1 < i \leq j < n$.}
    \end{cases}
\]
Therefore, both $a_1$ and $a_n$ should appear in the product
$\lambda(\pm a_1, a_i, \ldots, a_j, \pm a_n)$,
but just once each. If $n = 1$ then $a_1 = a_n$ should also appear just once.

Of course, there is an analogue of Theorem \thmref{rightproduct}.

\begin{theorem}
\thmlabel{leftproduct}
Let $a_1,\ldots,a_n \quad (n > 0)$
be a sequence of elements of the right product right loop
$S$, such that $a_{p_1}, \ldots , a_{p_m} \quad (1 \leq p_1 < \ldots < p_m \leq n;
\quad 0 \leq m \leq n)$ and only them among $a_1, \ldots, a_n$ are nonidempotents.
Then $\lambda(a_1,\ldots,a_n) = \lambda(\pm a_1, a_{p_1},\ldots,a_{p_m},\pm a_n)$.
\end{theorem}

\begin{proof}
The proof is by induction on $n$.

\noindent (1) $n = 1$.

\noindent If $m = 0$ then $\lambda(a_1) = a_1 = \lambda(\pm a_1, \pm a_1)$.

\noindent If $m = 1$ then $\lambda(a_1) = a_1 = \lambda(\pm a_1, a_1, \pm a_1)$.

\noindent (2) $n > 1$.

\noindent (2a) Let $1 = p_1 \, , p_{m-1} = n-1 \, , p_m = n$
(i.e. $a_1, a_{n-1}, a_n \notin E_S$). Then, using induction argument,
$\lambda(a_1,\ldots,a_n)
= \lambda(a_1,\ldots,a_{n-1}) \cdot a_n
= \lambda(\pm a_1,a_{p_1}, \ldots, a_{p_{m-1}}, \pm a_{n-1}) \cdot a_{p_m}
= \lambda(a_{p_1},\ldots,a_{p_{m-1}}) \cdot a_{p_m}
= \lambda(a_{p_1}, \ldots, a_{p_{m-1}}, a_{p_m})
= \lambda(\pm a_1, a_{p_1}, \ldots, a_{p_m}, \pm a_n)$.

\noindent (2b) Let $1 = p_1 \, , p_m = n-1$ (i.e.
$a_1, a_{n-1} \notin E_S \, ;a_n \in E_S$).
Then, using induction argument again,
we get
$\lambda(a_1,\ldots,a_n)
= \lambda(a_1,\ldots,a_{n-1}) \cdot a_n
= \lambda(\pm a_1, a_{p_1}, \ldots, a_{p_m}, \pm a_{n-1}) \cdot a_n
= \lambda(a_{p_1}, \ldots, a_{p_m}) \cdot a_n
= \lambda(a_{p_1}, \ldots, a_{p_m}, a_n)
= \lambda(\pm a_1, a_{p_1}, \ldots, a_{p_m}, \pm a_n)$.

\noindent (2c) Let $1 = p_1 \, , p_{m-1} < n-1 \, , p_m = n$
(i.e. $a_1, a_n \notin E_S \, ;a_{n-1} \in E_S$).
Then, by the induction argument and (RL),
$\lambda(a_1,\ldots,a_n)
= \lambda(a_1,\ldots,a_{n-1}) \cdot a_n
= \lambda(\pm a_1,a_{p_1}, \ldots, a_{p_{m-1}}, \pm a_{n-1}) \cdot a_n
= \lambda(a_{p_1},\ldots,a_{p_{m-1}}, a_{n-1}) \cdot a_{p_m}
= \lambda(a_{p_1},\ldots,a_{p_{m-1}}) a_{n-1} \cdot a_{p_m}
= \lambda(a_{p_1},\ldots,a_{p_{m-1}}) \cdot a_{p_m}
= \lambda(a_{p_1}, \ldots, a_{p_{m-1}}, a_{p_m})
= \lambda(\pm a_1, a_{p_1}, \ldots, a_{p_m}, \pm a_n)$.

\noindent (2d) Let $1 = p_1 \, , p_{m-1} < n$
(i.e. $a_1 \notin E_S \, ;a_{n-1}, a_n \in E_S$).
Then,
$\lambda(a_1,\ldots,a_n)
= \lambda(a_1,\ldots,a_{n-1}) \cdot a_n
= \lambda(\pm a_1, a_{p_1}, \ldots, a_{p_m}, \pm a_{n-1}) \cdot a_n
= \lambda(a_{p_1}, \ldots, a_{p_m}, a_{n-1}) \cdot a_n
= \lambda(a_{p_1}, \ldots, a_{p_m}) a_{n-1} \cdot a_n
= \lambda(a_{p_1}, \ldots, a_{p_m}) \cdot a_n
= \lambda(a_{p_1}, \ldots, a_{p_m}, a_n)
= \lambda(\pm a_1, a_{p_1}, \ldots, a_{p_m}, \pm a_n)$.

The remaining cases of (2) in which $p_1 \neq 1$, \emph{i.e.}, $a_1 \in E_S$
can be proved analogously.
\end{proof}

In right product $\{$left, right$\}$ loops, Theorems \thmref{rightproduct}
and \thmref{leftproduct} give us the means to reduce $\{$right, left$\}$
products. The result is much stronger in right product loops.

\begin{definition}
\deflabel{ExternalNeutral}
Let $(S; \cdot, \ldiv, \rdiv)$ be a right product quasigroup and $1 \notin S$.
By $S^1$ we denote a triple magma with operations extending
$\cdot, \ldiv, \rdiv$ to $S \cup \{ 1 \}$ in the following way:
$x \circ y \quad (\circ \in \{ \cdot, \ldiv, \rdiv \})$ remains as before
if $x, y \in S$. If $x = 1$ then $x \circ y = y$ and if $y = 1$ then
$x \circ y = x$.
\end{definition}

Note that the new, extended operations $\cdot, \ldiv, \rdiv$ are well defined
and that $1$ is the neutral element for all three.

\begin{lemma}
\lemlabel{id3ass}
Let $a_1, \ldots, a_n \quad (n > 0)$ be a sequence of elements
of a right product loop $S$ such that $a_n$ is an idempotent and
$p(a_1, \ldots, a_n)$ some product of $a_1, \ldots, a_n$ (in that order)
with an arbitrary (albeit fixed) distribution of brackets . Then
$p(a_1, \ldots, a_n) = p(a_1, \ldots, a_{n-1},1) \cdot a_n$.
\end{lemma}

\begin{proof}
First, note that if $e$ is an idempotent then
$x \cdot ye = xy \cdot e$ for all $x, y \in S$. Namely,
if $e \in E_S$ then there is a $z \in S$ such that $e = z \rdiv z$
(for example $z = e$ is one). The identity
$x \cdot y(z \rdiv z) = xy \cdot (z \rdiv z)$ is true in all right product
loops as it is true in all loops and all right zero semigroups.

The proof of the lemma is by induction on $n$.

\noindent (1) $ n = 1$.

\noindent $a_1 = a_n$ is an idempotent, so $p(a_1) = a_1 = 1 \cdot a_1 = p(1) \cdot a_1$.

\noindent (2) $n > 1$.

\noindent Let $p(a_1, \ldots, a_n) = q(a_1, \ldots, a_k) \cdot r(a_{k+1}, \ldots, a_n)$
for some $k \quad (1 \leq k \leq n)$.
By the induction hypothesis
$r(a_{k+1}, \ldots, a_n) = r(a_{k+1}, \ldots, a_{n-1}, 1) \cdot a_n$.
So $p(a_1, \ldots, a_n)
= q(a_1, \ldots, a_k) \cdot (r(a_{k+1}, \ldots, a_{n-1}, 1) \cdot a_n)
= (q(a_1, \ldots, a_k) \cdot r(a_{k+1}, \ldots, a_{n-1}, 1)) \cdot a_n
= p(a_1, \ldots, a_{n-1}, 1) \cdot a_n$.
\end{proof}

The following result is an improvement
of Theorems \thmref{rightproduct} and \thmref{leftproduct}.

\begin{theorem}
\thmlabel{product}
Let $a_1, \ldots, a_n$ and
$b_1, \ldots, b_n \quad (n > 0)$ be two sequences
of elements
of the right product loop $S$ (with some of $b_k$ possibly being $1$) such that
\[
b_k =
\begin{cases}
1 ; & \text{if $k < n$ and $a_k \in E_S$} \\
a_k ; & \text{if $k = n$ or $a_k \notin E_S$}
\end{cases}
\]
and let $p(a_1, \ldots, a_n)$ be as in Lemma \lemref{id3ass}. Then
$p(a_1, \ldots, a_n) = p(b_1, \ldots, b_n)$.
\end{theorem}

\begin{proof}
The proof of the Theorem is by induction on $n$.

\noindent (1) $n = 1$.

\noindent There is only one product $p(a_1) = a_1$ and, irrespectively
of whether $a_1$ is idempotent or not, $b_1 = a_1$.
Therefore $p(a_1) = p(b_1)$.

\noindent (2) $n > 1$.

\noindent  Let $p(a_1, \ldots, a_n) = q(a_1, \ldots, a_k) \cdot r(a_{k+1}, \ldots, a_n)$
for some $k \quad (1 \leq k \leq n)$.
By the induction hypothesis
$q(a_1, \ldots, a_k) = q(b_1, \ldots, b_{k-1}, a_k)$ and
$r(a_{k+1}, \ldots, a_n) = r(b_{k+1}, \ldots, b_n)$.

If $a_k$ is nonidempotent then $a_k = b_k$ and
$p(a_1, \ldots, a_n) = q(b_1, \ldots, b_k) \cdot r(b_{k+1}, \ldots, b_n)
= p(b_1, \ldots, b_n)$.

If $a_k$ is idempotent then $b_k = 1$ and by the
Lemma \lemref{id3ass}
$p(a_1, \ldots, a_n)
= q(b_1, \ldots, b_{k-1}, a_k) \cdot r(b_{k+1}, \ldots, b_n)
= (q(b_1, \ldots, a_{k-1}, 1) \cdot a_k) \cdot r(b_{k+1}, \ldots, b_n)
= q(b_1, \ldots, b_k) \cdot r(b_{k+1}, \ldots, b_n)
= p(b_1, \ldots, b_n)$.
\end{proof}

The following corollary is an analogue of (\cite{KrapezRL}, Theorem 2.4).

\begin{corollary}
\corlabel{product2}
Let $a_1,\ldots,a_n$ be a sequence of elements of the
right product loop $S$, such that at most two of them are
nonidempotents.  Then all products of $a_1,\ldots,a_n$, in that
order, are equal to the following product of at most three of
them: First -- nonidempotents of $a_1,\ldots,a_{n-1}$ if any (the
one with the smaller index first) and then $a_n$ if it is not used already.
\end{corollary}

In right product pointed loops we need not use $1$.

\begin{theorem}
\thmlabel{product1}
Let $a_1, \ldots, a_n$ and
$b_1, \ldots, b_n \quad (n > 0)$ be two sequences
of elements
of the right product pointed loop $S$ with the distinguished element $e$
such that
\[
b_k =
\begin{cases}
e ; & \text{if $k < n$ and $a_k \in E_S$} \\
a_k ; & \text{if $k = n$ or $a_k \notin E_S$}
\end{cases}
\]
and let $p(a_1, \ldots, a_n)$ be some product of $a_1, \ldots, a_n$.
Then
$p(a_1, \ldots, a_n) = p(b_1, \ldots, b_n)$.
\end{theorem}

\section{Independence of axioms}
\seclabel{indep}

Finally, we consider the independence of the axioms (\eqA)
for right product quasigroups.

It is well-known that the quasigroup axioms (\eqQ{1})--(\eqQ{4})
are independent. It follows that axioms (\eqA{1}) and (\eqA{2})
are independent. To give just one concrete example, here is a
model in which (\eqQ{2}) $=$ (\eqA{2}) fails.

\begin{example}
\exmlabel{notA2}
The model $(\mathbb{Z}; \cdot, \ldiv, \rdiv)$ where $x\cdot y = x + y$,
$x \rdiv y = x - y$ and $x \ldiv y = \max\{y-x,0\}$ is a
left quasigroup satisfying (\eqQ{1}) but not (\eqQ{2}),
and hence satisfies (\eqA{1}), (\eqA{3}),
(\eqA{4}) and (\eqA{5}), but not (\eqA{2}).
\end{example}

As it turns out, the independence of the remaining axioms
can be easily shown by models of size $2$. These were found using
\textsc{Mace4} \cite{McCune}.

\begin{example}
\exmlabel{notA3}
Table \tabref{notA3} is a model satisfying
(\eqA{1}), (\eqA{2}), (\eqA{4}), (\eqA{5}), but not (\eqA{3}).
\begin{table}[htb] \centering
\[
\begin{array}{r|rr}
\cdot & 0 & 1\\
\hline
    0 & 0 & 1 \\
    1 & 0 & 1
\end{array}
\qquad
\begin{array}{r|rr}
\ldiv & 0 & 1\\
\hline
    0 & 0 & 1 \\
    1 & 0 & 1
\end{array}
\qquad
\begin{array}{r|rr}
\rdiv  & 0 & 1\\
\hline
    0 & 1 & 0 \\
    1 & 1 & 0
\end{array}
\]
\caption{(\eqA{1}), (\eqA{2}), (\eqA{4}), (\eqA{5}), but not (\eqA{3})}
\tablabel{notA3}
\end{table}
\end{example}

\begin{example}
\exmlabel{notA4}
Table \tabref{notA4} is a model satisfying
(\eqA{1}), (\eqA{2}), (\eqA{3}), (\eqA{5}), but not (\eqA{4}).
\begin{table}[htb]  \centering
\[
\begin{array}{r|rr}
\cdot & 0 & 1\\
\hline
    0 & 0 & 1 \\
    1 & 1 & 0
\end{array} \qquad
\begin{array}{r|rr}
\ldiv & 0 & 1\\
\hline
    0 & 0 & 1 \\
    1 & 1 & 0
\end{array} \qquad
\begin{array}{r|rr}
\rdiv & 0 & 1\\
\hline
    0 & 1 & 0 \\
    1 & 0 & 1
\end{array}
\]
\caption{(\eqA{1}), (\eqA{2}), (\eqA{3}), (\eqA{5}), but not (\eqA{4})}
\end{table}
\tablabel{notA4}
\end{example}

\begin{example}
\exmlabel{notA5}
Table \tabref{notA5} is a model satisfying
(\eqA{1}), (\eqA{2}), (\eqA{3}), (\eqA{4}), but not (\eqA{5}).
\begin{table}[htb]  \centering
\[
\begin{array}{r|rr}
\cdot & 0 & 1\\
\hline
    0 & 1 & 0 \\
    1 & 1 & 0
\end{array} \qquad
\begin{array}{r|rr}
\ldiv & 0 & 1\\
\hline
    0 & 1 & 0 \\
    1 & 1 & 0
\end{array} \qquad
\begin{array}{r|rr}
\rdiv & 0 & 1\\
\hline
    0 & 1 & 0 \\
    1 & 1 & 0
\end{array}
\]
\caption{(\eqA{1}), (\eqA{2}), (\eqA{3}), (\eqA{4}), but not (\eqA{5})}
\tablabel{notA5}
\end{table}
\end{example}

%%%%%%%%%%%%%%%%%%%%%%%%%%%%%%%

\end{document}